\documentclass[a4paper,12pt]{article}

\usepackage{amsmath}
\usepackage{amsthm}
\usepackage{amssymb}  
\usepackage{latexsym} 
\usepackage[all]{xy}
\usepackage{comment}
\usepackage{rotating}
\usepackage{slashbox}
\usepackage{multirow}
\usepackage{rotating}

\newcommand{\Char}{\operatorname{char}}

\newcommand{\Mathematica}{{\sf Mathematica }}
\newcommand{\Sage}{{\sf Sage }}

\newcommand{\isom}{ \cong }

\newcommand{\PP}{{\mathbb P}}
\newcommand{\Q}{{\mathbb Q}}

\newcommand{\rank}{\operatorname{rank}}

\newcommand{\tor}{\operatorname{tor}}

\newcommand{\Z}{{\mathbb Z}}

\newfont{\wncyr}{wncyr10 at 12pt}
\newfont{\wncyrten}{wncyr10 at 10pt}

\newenvironment{Proof}{\par\noindent{\sc Proof:}}%
                      {\hspace*{\fill}\nobreak$\Box$\par\medskip}
                       {\hspace*{\fill}\nobreak$\Box$\par\medskip}
\newenvironment{myitemize}
{\begin{itemize}
\setlength{\itemsep}{1pt}
\setlength{\parskip}{0pt}
\setlength{\parsep}{0pt}}
{\end{itemize}}
\newenvironment{myenumerate}
{\begin{enumerate}
\setlength{\itemsep}{1pt}
\setlength{\parskip}{0pt}
\setlength{\parsep}{0pt}}
{\end{enumerate}}

\newtheorem{Proposition}{Proposition}[section]
\newtheorem{Theorem}[Proposition]{Theorem}
\newtheorem{Lemma}[Proposition]{Lemma}
\newtheorem{Corollary}[Proposition]{Corollary}

\theoremstyle{definition}

\newtheorem{Remark}[Proposition]{Remark}

\newcounter{nootje}
\begin{document}
\normalsize
\title{Periodic continued fractions and elliptic curves over quadratic fields}
\author{Mohammad Sadek}
\date{}
\maketitle
\let\thefootnote\relax\footnote{Mathematics Subject Classification: 11A55, 11J70}
\begin{abstract}{\footnotesize Let $f(x)$ be a square free quartic polynomial defined over a quadratic field $K$ such that its leading coefficient is a square. If the continued fraction expansion of $\displaystyle \sqrt{f(x)}$ is periodic, then its period $n$ lies in the set \[\{1,2,3,4,5,6,7,8,9,10,11,12,13,14,15,17,18,22,26,30,34\}.\] We write explicitly all such polynomials for which the period $n$ occurs over $K$ but not over $\Q$ and $\displaystyle n\not\in\{ 13,15,17\}$. Moreover we give necessary and sufficient conditions for the existence of such continued fraction expansions with period $13, 15$ or $17$ over $K$.}
\end{abstract}

\section{Introduction}

Let $E$ be an elliptic curve defined over a field $K$ whose characteristic is different from $2,3$. One can describe $E$ using an equation of the form $y^2=f(x)$ where $f(x)$ is a square free quartic polynomial whose leading coefficient is a square in $K$. In \cite{Razar} the authors were able to prove that the continued fraction expansion of $\sqrt{f(x)}$ is periodic if and only if the point $\infty^+-\infty^-$ is of finite order on $E$. Furthermore the period of the continued fraction expansion can be determined once the order of the point $\infty^+-\infty^-$ is known.

Given a square free quartic polynomial $f(x)=a_0x^4+a_1x^3+a_2x^2+a_3x+a_4$ defined over $K$ whose leading coefficient $a_0$ is a square, the curve defined by the equation $C:y^2=f(x)$ has a rational point, namely $(1:\pm\sqrt{a_0}:0)$. Therefore the curve $C$ is isomorphic over $K$ to its Jacobian elliptic curve $E$. The continued fraction expansion of $\sqrt{f(x)}$ is periodic if and only if the order of $\infty^+-\infty^-$ in $E(K)$ is finite. In fact if the order of $\infty^+-\infty^-$ is $n$ then the period of the continued fraction is either $n-1$ or $2(n-1)$ where the second case occurs only if $n$ is even.

The above argument leads one to study elliptic curves with torsion points in order to investigate quartic polynomials $f(x)$ where the continued fraction expansion of $\sqrt{f(x)}$ is periodic. An elliptic curve with a $K$-rational torsion point of order $n$ can be written in Tate's normal form; namely
\[E_{b,c}:y^2+(1-c)xy-by=x^3-bx^2,\;\textrm{ where }b,c\in K.\] If $f(x)$ has a square leading coefficient such that the continued fraction expansion of $\sqrt{f(x)}$ is periodic then there exist $b,c\in K$ such that the curve $C:y^2=f(x)$ is isomorphic to $E_{b,c}$.

In \cite{PeriodicCF} Alfred van der Poorten wrote explicitly all square free quartic polynomials $f(x)$ over $\Q$ with a square leading coefficient such that $\sqrt{f(x)}$ is periodic. Following Mazur's classification of torsion points of elliptic curves over $\Q$ the possible periods are
\[\{1,2,3,4,5,6,7,8,9,10,11,14,18,22\}.\]  In fact it was shown that all of these periods occur over $\Q$ except for $9$ and $11$ as there is no such polynomial over $\Q$ such that the continued fraction expansion of $\sqrt{f(x)}$ is of period $9$ nor $11$.

In this article we write down all square free quartic polynomials $f(x)$ with a square leading coefficient such that the continued fraction expansion of $\sqrt{f(x)}$ is periodic over some quadratic field $K$. According to the classification of torsion points of elliptic curves over quadratic fields the possible periods are the ones over $\Q$ together with  \[\{9,11,12,13,15,17,26,30,34\}.\] We prove that the periods $9,11$ occur over some quadratic fields. Moreover we display all quartic polynomials that give rise to the periods $12,26,30,34$. In addition we present the quadratic fields with the smallest absolute value of their discriminants over which these periods occur. Finally we give necessary and sufficient conditions for the odd periods $13, 15, 17$ to occur over a quadratic field $K$. More precisely we show that the period $n$, $n=13, 15, 17$, occurs if and only if there exists a $z\in K$ such that $z^2=\alpha(t,s)$ for some $K$-rational point $(t,s)$ lying on the modular curve $X_1(n+1)$.

\section{Continued fractions and elliptic curves}

Let $E$ be an elliptic curve defined over a field $K$ with $\Char K\ne 2$. Let $\pi:E\to \PP^1$ be a degree $2$ separable map. Let $x\in K(\PP^1)$ have a simple pole at a point $\infty\in \PP^1$. Thus one can identify $K(x)$ with $K(\PP^1)$. Choosing $y\in K(E)$ with a pole of order $2$ at $\infty^+$, one can construct a hyperelliptic equation $y^2=f(x)$ describing $E$ where $f(x)\in K[x]$ is a polynomial of degree $4$. We can assume moreover that the leading coefficient of $f(x)$ is a square. The pullback of $\infty\in K(\PP^1)$ under $\pi$ consists of two unramified points $\infty^+$ and $\infty^-$ defined over $K$. Moreover $\infty^+-\infty^-$ is a point on the Jacobian curve of $E$ which is $E$ itself. Now $1/x$ is a uniformising parameter at $\infty^+$, and defines an embedding $K(E)\hookrightarrow L:=K((1/x))$.

Let $\alpha$ be an element in $K(E)$. Then the image of $\alpha$ in $L$ is the Laurent series of $\alpha$. The element $\alpha$ has a continued fraction expansion
\begin{equation}
  \alpha(x) = a_0(x) + \cfrac{1}{a_1(x)
          + \cfrac{1}{a_2(x)
          + \cfrac{1}{a_3(x) + {}_{\ddots} } } }
\end{equation}
which can be written shortly as $[a_0,a_1,a_2,a_3,\ldots]$. The continued fraction algorithm is defined as follows: $a_0$ is the polynomial part of $\alpha$. The {\em partial quotient} $a_h$ is defined to be the polynomial part of $\alpha_h$, then the $(h+1)$-th {\em complete quotient} is $\alpha_{h+1}=(\alpha_h-a_h)^{-1}$ which is again an element in $L$. The $h$-th convergent is defined to be \[\frac{p_h}{q_h}=[a_0,a_1,\ldots,a_h],\;h=0,1,2,\ldots.\]

The complete quotient $\alpha_h$ of $\alpha$ is of the form $\displaystyle\alpha_h=\left(P_h+\sqrt{f(x)}\right)/Q_h$ where $P_{h+1}+P_h=a_{h}Q_h$, and $Q_{h+1}Q_h=f(x)-P_{h+1}^2.$ Furthermore $\deg P_{h+1}=2$ and $\deg Q_h\le 1$ for all $h=0,1,2,\ldots$.
The partial quotients $a_h$ of $\alpha$, $h\ge 1,$ are polynomials of degree $1$. When $B,C$ are nonzero in $K[x]$ one obtains the following identity
\[B[Ca_0,Ba_1,Ca_2,Ba_3,Ca_4,\ldots]=C[Ba_0,Ca_1,Ba_2,Ca_3,Ba_4,\ldots],\]
see \cite[p. 140]{Schmidt}.

A continued fraction of some $\alpha\in K(E)$ is said to be {\em periodic} if there is a positive integer $n$ such that $a_{h+n}=a_h$ for every $h>1$. The smallest such $n$ is called the {\em period} of the continued fraction.

\begin{Lemma}
\label{lem:periodic}
Let $y^2=f(x)$ be a hyperelliptic equation describing $E$ where $f(x)\in K[x]$ is a square free polynomial of degree $4$ whose leading coefficient is a square. The following statements are equivalent.
\begin{myenumerate}
\item The point $\infty^+-\infty^-$ is a torsion point on $E$.
\item $y=\sqrt{f(x)}$ has a periodic continued fraction expansion.
\end{myenumerate}
\end{Lemma}
\begin{Proof}
This is the main result of \cite{Razar}, see \cite[Lemma 7]{Schmidt}.
\end{Proof}
The period of the continued fraction expansion of $y$ is determined by the following proposition, see \cite[Corollary 4.4]{Razar}.

\begin{Proposition}
\label{prop:periodlength}
If the order of $\infty^+-\infty^-$ is $n$ then the period of the continued fraction of $y$ is either $n-1$ or $2(n-1)$. Moreover the second case holds only if $n$ is even.
\end{Proposition}


The following lemma describes the partial quotients of the continued fraction in Lemma \ref{lem:periodic}.

 \begin{Lemma}
\label{lem:periodicexpansion}
Let $f(x)$ be a quartic polynomial over $K$ whose leading coefficient is a square. Assume moreover that $f(x)$ has no repeated roots. Assume that the continued fraction expansion of $\displaystyle \sqrt{f(x)}$ is periodic. Then the continued fraction expansion of $\displaystyle \sqrt{f(x)}$ is given by:
  \[\sqrt{f(x)}=\left[ a_0, \overline{a_1,\ldots,a_{r-1},2a_0/k,ka_1,\ldots,a_{r-1}/k,2a_0}\right]\]
  where $a_0$ is of degree $2$, and $a_i$ is of degree $1$, $1\le i\le r-1$. Furthermore $a_i=k^{m_i}a_{r-i}$, where $m_i=(-1)^{i+1}$, $1\le i< (r+1)/2$. Hence if $r$ is even then $k=1$.
\end{Lemma}
\begin{Proof}
This is \cite[Proposition 3]{ellipticfunctionfields}.
\end{Proof}

\begin{Remark}
\label{rem1}
In Lemma \ref{lem:periodicexpansion} if $k=1$, then the period length is $r$. If $k\ne 1$, in particular $r$ is odd, then $\mu y$ has period length $r$ if and only if $\mu^2=1/k$. This follows because the continued fraction of $\mu y$ is given by \[\mu y=\left[ \mu a_0, \overline{a_1/\mu,\mu a_2,\ldots,\mu a_{r-1},2a_0/(k\mu),k\mu a_1,\ldots,a_{r-1}/(k\mu),2a_0\mu}\right],\] hence $\mu=1/(k\mu).$
\end{Remark}
\begin{Remark}
\label{rem:partial quotients}
The explicit description of the partial quotients of the continued fraction expansion of $\sqrt{f(x)}$, where $f(x)=(x^2+u)^2-4v(x+w)$, can be found in \cite[\S4]{PeriodicCF} and is given as follows:
\begin{eqnarray*}
a_0&=&x^2+u,\;\\
a_1&=& 2(x-w)/v,\;\\
a_i&=&2(x-c_i)/b_i,\;i\ge 2,
\end{eqnarray*}
where \begin{eqnarray*}b_{2h}=\frac{s_3s_5\ldots s_{2h-1}}{s_2s_4\ldots s_{2h}}\textrm{ and }\;b_{2h+1}=4v\frac{s_2s_4\ldots s_{2h}}{s_3s_5\ldots s_{2h+1}},
\end{eqnarray*}
\begin{eqnarray*}
s_2&=&1;\\
s_3&=&v/(-2w+1);\\
s_h&=&v/(s_{h-1}(s_{h-1}-1)s_{h-2}),\;h\ge 4.
\end{eqnarray*}
If the continued fraction of $\sqrt{f(x)}$ is periodic with period $n$ then $c_h=c_{n-h}, 1\le h \le n-1$.
\end{Remark}

\section{Torsion points of elliptic curves}

Since Mazur's classification of torsion points of elliptic curves defined over $\Q$, an enormous amount of research has been directed toward studying torsion points of elliptic curves over number fields.
\begin{Theorem}
\label{mazur}
Let $E$ be an elliptic curve defined over $\Q$. Then the torsion subgroup $E_{\tor}(\Q)$ is isomorphic to one of the following 15 groups:
 \begin{eqnarray*}
 \Z/n\Z,\;1\le n\le 12,\;n\ne11;\textrm{ or }
 \Z/2\Z\times\Z/2n\Z,\;1\le n\le 4.
 \end{eqnarray*}
\end{Theorem}
The theorem above together with Proposition \ref{prop:periodlength} implies that if $y^2=f(x)$ is a quartic model for $E$ such that $\infty^+-\infty^-$ has finite order then $\sqrt{f(x)}$ has continued fraction expansion with period $r$ where $r\in \{1,2,3,4,5,6,7,8,9,10,11,14,18,22\}$. In \cite{ellipticfunctionfields} it was proved that the period $r$ does not take the values $9,11$.

For quadratic extensions of $\Q$, a series of papers by Kamienny, Kenku, and Momose leads to the following analogue of Mazur's classification of torsion points of elliptic curves.

\begin{Theorem}
\label{thm:quadratictorsion}
Let $K$ be a quadratic extension of $\Q$. Let $E$ be an elliptic curve defined over $K$. Then the torsion subgroup $E_{\tor}(K)$ is isomorphic to one of the following 26 groups:
 \begin{eqnarray*}
 \Z/n\Z&,&\;1\le n\le 18,\;n\ne17;\\
 \Z/2\Z\times\Z/2n\Z&,&\;1\le n\le 6;\\
 \Z/3\Z\times\Z/3m\Z&,&\;1\le n\le 2;\\
 \Z/4\Z\times\Z/4\Z&.&
 \end{eqnarray*}
\end{Theorem}
\begin{Remark}
\label{rem2}
The theorem above together with Proposition \ref{prop:periodlength} indicate that if $y^2=f(x)$ is a quartic model for $E/K$ then if $\infty^+-\infty^-$ is a point of finite order on $E$ then the continued fraction of $\sqrt{f(x)}$ is periodic, and the period is in the set \[S=\{1,2,3,4,5,6,7,8,9,10,11,12,13,14,15,17,18,22,26,30,34\}.\] It follows that the periods that may occur over $K$ but not over $\Q$ are in the set \[S'=\{9,11,12,13,15,17,26,30,34\}.\]
\end{Remark}
In this paper we investigate whether all the values in the set $S'$ are realized over some quadratic field. In other words, does every $r\in S'$ appear as the period of the continued fraction of $\sqrt{f(x)}$ for some square free polynomial $f(x)\in K[x]$ of degree $4$ whose leading coefficient is a square, where $K$ is a quadratic field?

Tate's normal form is a Weierstrass equation describing an elliptic curve with non-trivial torsion points, namely it is given by
\begin{eqnarray*}
y^2+(1-c)xy-by=x^3-bx^2
\end{eqnarray*}
where the point $(0,0)$ is a point of finite maximal order which is neither $2$ nor $3$. Kubert was the first to list all elliptic curves over $\Q$ with a torsion point in Tate's normal form, see \cite{Kubert}. More precisely, a parametrization of the coefficients $b,c$ were given explicitly in terms of a parameter $t$ for every elliptic curve in the list of Theorem \ref{mazur}. These normal forms were used in \cite{ellipticfunctionfields} to produce all periodic continued fractions of $\sqrt{f(x)}$, where $f(x)$ is a square free polynomial of degree 4 over $\Q$.

In order to extend the work of \cite{ellipticfunctionfields} over quadratic fields, one needs to seek a parametrization for elliptic curves with torsion points over quadratic fields. This can be found in \cite{Rabarison}. Let $Y_1(M,N)$ be the affine modular curve parametrizing isomorphism classes of triples $(E,P_M,P_N)$ where $P_M,P_N$ are points on the elliptic curve $E$ such that $\langle P_M\rangle\times\langle P_N\rangle\isom \Z/M\Z\times\Z/N\Z$ with $M|N$. Let $X_1(M,N)$ be the compactification of $Y_1(M,N)$, namely $X_1(M,N)$ is obtained from $Y_1(M,N)$ by adding the cusps. For simplicity we are going to write $X_1(N)$ instead of $X_1(1,N)$.

The following theorem is \cite[Theorem 2]{NajmanNumberTheorey}. It provides the complete list of orders of torsion points of elliptic curves defined over either $\Q(\sqrt{-1})$ or $\Q(\sqrt{-3})$.

\begin{Theorem}
\label{thm:Qi}
\begin{myenumerate}
\item Let $E$ be an elliptic curve over $\Q(\sqrt{-1})$. Then $E\left(\Q(\sqrt{-1})\right)_{tor}$ is either one of the groups from Mazur's Theorem or $\Z/4\Z\times\Z/4\Z$.
\item  Let $E$ be an elliptic curve over $\Q(\sqrt{-3})$. Then $E\left(\Q(\sqrt{-3})\right)_{tor}$ is either one of the groups from Mazur's Theorem, $\Z/3\Z\times\Z/3\Z$ or  $\Z/3\Z\times\Z/6\Z$.
\end{myenumerate}
\end{Theorem}

 In fact, the only quadratic field over which torsion $\Z/4\Z\times\Z/4\Z$ occurs is $\Q(\sqrt{-1})$, and the only quadratic field over which torsion $\Z/3\Z\times\Z/3\Z$ and $\Z/3\Z\times\Z/6\Z$ occur is $\Q\left(\sqrt{-3}\right)$.

For the convenience of the reader we write down the defining equations of the curves $X_1(M,N)$ that correspond to elliptic curves with torsion points over a quadratic field, but not all of them are defined over the rational field. We use the following defining equations that were collected from different sources in \cite{Najman}.
{\center
{\scriptsize
\begin{tabular}{|c||c|c|}
               \hline
               Curve & Equation & Cusp \\
               \hline
               $X_1(11)$ & $y^2-y=x^3-x$ & $x(x - 1)(x^5 - 18x^4 + 35x^3 - 16x^2 - 2x + 1) = 0$ \\
               \hline
               $X_1(13)$ & $y^2 = x^6 - 2x^5 + x^4 - 2x^3 + 6x^2 - 4x + 1$ & $x(x - 1)(x^3 - 4x^2 + x + 1) = 0$ \\
               \hline
               $X_1(14)$ & $y^2 + xy + y = x^3 - x$ & $x(x - 1)(x + 1)(x^3 - 9x^2 - x + 1)(x^3 - 2x^2 - x + 1) = 0$ \\
               \hline
               $X_1(15)$ & $y^2 + xy + y = x^3 + x^2$ & $x(x + 1)(x^4 + 3x^3 + 4x^2 + 2x + 1)(x^4 - 7x^3- 6x^2 + 2x + 1) = 0$ \\
               \hline
               $X_1(16)$ & $y^2 = x(x^2 + 1)(x^2 + 2x - 1)$ & $x(x - 1)(x + 1)(x^2 - 2x - 1)(x^2 + 2x - 1) = 0$ \\
               \hline
               $X_1(18)$ & $y^2 = x^6 + 2x^5 + 5x^4 + 10x^3 + 10x^2 + 4x + 1$& $x(x + 1)(x^2 + x + 1)(x^2 - 3x - 1) = 0$ \\
               \hline
               $X_1(2,10)$ & $y^2 = x^3 + x^2 - x$ & $x(x - 1)(x + 1)(x^2 + x - 1)(x^2 - 4x - 1) = 0$ \\
               \hline
               $X_1(2,12)$ & $y^2 = x^3 - x^2 + x$ & $x(x - 1)(2x - 1)(2x^2 - x + 1)(3x^2 - 3x - 1)(6x^2 - 6x - 1) = 0$ \\
               \hline
             \end{tabular}}}
             \vskip10pt
In order to find elliptic curves over $K$ with torsion $\Z/n\Z\times\Z/m\Z$, one checks whether $X_1(M,N)$ has a $K$-rational point which is not a cusp.

\section{Quartic models of elliptic curves}

We write down the isomorphisms between Tate's elliptic curve $E_{b,c}:y^2+(1-c)xy-by=x^3-bx^2$ and an elliptic curve in short Weierstrass normal form. These can be found for example in \cite{ellipticfunctionfields}. This is performed via the following change of variables $x=u^2x'+r,y=u^3y'+su^2x'+t$, where

{\begin{center}
$\begin{array}{ll}
u=1,&r = -(c^2 - 2c - 4b + 1)/12,\\
s = (c - 1)/2,& t = -(c^3 - 3c^2 - (4b - 3)c - (8b + 1))/24.
\end{array}$
\end{center}}

 Now the elliptic curve $E_{b,c}$ is isomorphic to
  \[E:y^2=x^3+Ax+B,A,B\in K.\]
  Moreover the torsion point $(0,0)$ is transformed to a torsion point $P =(c^2 -2c - 4b + 1)/12, −b/2)$ of maximal order. In addition any other elliptic curve described by a short Weierstrass equation and isomorphic to $E_{b,c}$ will be isomorphic to $E$ via a transformation of the form $x=u^2x',y=u^3y'$.

 Let $v^2=u^3+Au+B,\;A,B\in K,$ be a Weierstrass equation describing an elliptic curve $E$ with a rational point $P(a,b)\in E(K)$. One can write a quartic model describing $E$ for which the points $P$ and $O$ are the points at infinity using the following transformation that can be found in \cite{Razar}. Namely, one has
 \begin{eqnarray}
 \label{eq2}
 (u,v)&\mapsto&(x,y)=\left(\frac{v+b}{u-a},\;2u+a-\left((v+b)/(u-a)\right)^2\right);\nonumber\\
 (x,y)&\mapsto&(u,v)=\left((x^2+y-a)/2,\;(x^3+xy-3ax-2b)/2\right)
 \end{eqnarray}
define a birational equivalence between $E$ and the curve $E_P$ described by
\begin{eqnarray*}
y^2=x^4-6ax^2-8bx-4A-3a^2,\textrm{ where }B=b^2-a^3-Aa.\\
\end{eqnarray*}
Furthermore the point $P$ and the point at infinity $O$ on $E$ are sent to the points at infinity on $E_P$. We conclude with the following proposition.
\begin{Proposition}
\label{prop:parametrizing quartics}
Let $E$ be an elliptic curve over a field $K$ such that $\Char K\ne 2,3$. Assume moreover that $E$ has a $K$-rational point of finite order $m$. Then there exist $b,c\in K$ such that $E$ is $K$-isomorphic to an elliptic curve described by the following equation \[y^2=\left(x^2+u_m(b,c)\right)^2-4v_{m}(b,c)\left(x+w_m(b,c)\right)\] where
\begin{eqnarray*}
u_m(b,c)&=&-3a(b,c)=-(c^2 -2c - 4b + 1)/4;\\
v_m(b,c)&=&b;\\
w_m(b,c)&=&-(A(b,c)+3a(b,c)^2)/b;\\
A(b,c)&=&(-16 b^2 - (-1 + c)^4 + 8 b (-2 + c + c^2))/48;\\
B(b,c)&=& (-64 b^3 + (-1 + c)^6 - 12 b (-1 + c)^3 (2 + c) +
   24 b^2 (5 + 2 c + 2 c^2))/864\\&=& (b/2)^2-a(b,c)^3-A(b,c)a(b,c).
\end{eqnarray*}
\end{Proposition}
\begin{Proof}
This follows from the argument above observing that any elliptic curve $E$ with a point of finite order can be described by a Tate's normal form equation $y^2+(1-c)xy-by=x^3-bx^2$ where $(0,0) $ is the torsion point of maximal order. Then using a birational equivalence one can describe $E$ by the following quartic model
\begin{eqnarray*}
y^2&=&x^4-6a(b,c)x^2-4bx-4A(b,c)-3a(b,c)^2\\
&=& (x^2-3a(b,c))^2-4b\left(x+\frac{A(b,c)+3a(b,c)^2}{b}\right).
\end{eqnarray*}
Now two birationally equivalent projective curves are isomorphic.
\end{Proof}
\begin{Corollary}
\label{cor:parametrizationofquartics}
Let $f(x)$ be a square free quartic polynomial defined over $K, \;\Char K\ne2,3,$ whose leading coefficient is a square. If the continued fraction expansion of $\sqrt{f(x)}$ is periodic then the curve $y^2=f(x)$ is $K$-isomorphic to an elliptic curve defined by the following equation \[y^2=\left(x^2+u_m(b,c)\right)^2-4v_{m}(b,c)\left(x+w_m(b,c)\right)\] where $u_m(b,c),v_m(b,c),w_m(b,c)$ are defined as in Proposition \ref{prop:parametrizing quartics}.
\end{Corollary}
\begin{Proof}
Let $f(x)=a_0x^4+a_1x^3+a_2x^2+a_3x+a_4$. The curve $C:y^2=f(x)$ is an elliptic curve since $(1:\pm\sqrt{a_0}:0)\in C(K)$. It follows that $C$ is $K$-isomorphic to its Jacobian elliptic curve $E$. Moreover since the continued fraction of $\sqrt{f(x)}$ is periodic then $\infty^+-\infty^-$ is of finite order in $E(K)$. One concludes using Proposition \ref{prop:parametrizing quartics}.
\end{Proof}

\section{Continued fraction with periods 9 and 11}

In \cite{ellipticfunctionfields} it was proved that there is no quartic polynomial $f(x)$ over $\Q$ with square leading coefficient such that the continued fraction of $\sqrt{f(x)}$ is periodic of period 9 or 11. We give explicit examples of square-free polynomials of degree $4$ whose square root has a continued fraction expansion with periods $9$ and $11$ over quadratic fields. In fact we find the quadratic fields with the smallest absolute value of the discriminant $|\Delta|$ over which these polynomials exist, and the quadratic fields with the smallest $|\Delta|$ over which infinitely many such polynomials exist.

\begin{Theorem}
\label{thm:period9}
The quadratic field $\Q\left(\sqrt{-1}\right)$ is the quadratic field with the smallest $|\Delta|$ over which there is a square free quartic polynomial $f(x)$ for which the continued fraction of $\sqrt{f(x)}$ has period $9$, whereas $\Q\left(\sqrt{2}\right)$ is the quadratic field with the smallest $|\Delta|$ over which there are infinitely many such quartic polynomials.
\end{Theorem}
\begin{Proof}
According to Proposition \ref{prop:periodlength} if the point $\infty^+-\infty^-$ on the hyperelliptic curve $y^2=f(x)$ has order $10$ then the period of the continued fraction of $\sqrt{f(x)}$ is either $9$ or $18$.

The parametrization of an elliptic curve with a torsion point of order 10 is given by $E_{b,c}:y^2+(1-c)xy-by=x^3-bx^2$ where \[(b,c)=\left(\frac{t^3(t-1)(2t-1)}{(t^2-3t+1)^2},-\frac{t(t-1)(2t-1)}{t^2-3t+1}\right).\]
According to Proposition \ref{prop:parametrizing quartics} the elliptic curve $E_{b,c}$ can be described by a quartic model given by $y^2=(X^2+u_{10})^2-4v_{10}(X+w_{10})$ where
\begin{eqnarray*}
u_{10}(t) = \frac{-4t^6 -16t^5 + 8t^4 + 8t^3 -4t+ 1}{4(t^2 -3t+ 1)^2},\\
v_{10}(t) = -\frac{t^3(t-1)(2t-1)}{(t^2 -3t+ 1)^2},\; w_{10}(t) = \frac{2t^3 -2t^2 -2t+ 1}{2(t^2 -3t+ 1)}
\end{eqnarray*}
with $\displaystyle t(t-1)(2t-1)(t^2-3t+1)\ne 0$ and $k_{10}(t)= -4t(t-1)(t^2 -3t+1)$. Now the continued fraction expansion of $y$ is given by
{\scriptsize\begin{eqnarray*}
y=\Big[X^2 -\frac{4t^6-16t^5+8t^4+8t^3-4t+1}{4(t^2-3t+1)^2} , \frac{(t^2-3t+1)^2}{2t^3(t-1)(2t-1)}\left(X - \frac{2t^3-2t^2-2t+1}{2(t^2-3t+1)}\right)
,2\left(X + \frac{2t^3-4t^2+4t-1}{2(t^2-3t+1)}\right),\\
-\frac{t^2-3t+1}{2t(t-1)(2t-1)}\left(X - \frac{2t^3-6t^2+4t-1}{2(t^2-3t+1)}\right)
,-\frac{2(t^2-3t+1)^2}{t}\left(X + \frac{2t-1}{2(t^2-3t+1)}\right), \frac{1}{2t^2(t-1)}\left(X + \frac{2t-1}{2(t^2-3t+1)}\right)
,\\
\frac{2(t^2-3t+1)^2}{2t-1}\left(X - \frac{2t^3-6t^2+4t-1}{2(t^2-3t+1)}\right)
, -\frac{ 1}{2t(t-1)(t^2-3t+1)}\left(X + \frac{2t^3-4t^2+4t-1}{2(t^2-3t+1)}\right),\\
\frac{-2(t^2-3t+1)^3}{t^2(2t-1)}\left(X - \frac{2t^3-2t^2-2t+1}{2(t^2-3t+1)}\right)
,-\frac{1}{2t(t-1)(t^2-3t+1)}\left(X^2- \frac{4t^6-16t^5+8t^4+8t^3-4t+1}{4(t^2-3t+1)^2}\right),\ldots\Big].
\end{eqnarray*}}

According to Lemma \ref{lem:periodicexpansion} the continued fraction of $\mu y$ is periodic of period $9$ if and only if $a_9(\mu y)/a_0(\mu y)=2$ where $a_i(\mu y)$ is the $i$-th partial quotient of $\mu y$. Equivalently, when $k_{10}=1/\mu^2$, see Remark \ref{rem1}. In other words, one needs to investigate $C_{10}(K)$ where $C_{10}:y^2=k_{10}:=-4t(t-1)(t^2 -3t+1),\;y\ne 0$. The Jacobian elliptic curve $E_{10}$ of $C_{10}$ is given by $y^2=x^3-x^2-x$.
One has $E_{10}(\Q(\sqrt{-1}))=\Z/6\Z$ which means that $E_{10}$ contains points which are not cusps over $\Q(\sqrt{-1})$, whereas $E_{10}(K)=\Z/2\Z$ when $K=\Q(\sqrt{-3})$. Therefore $\Q(\sqrt{-1})$ is the quadratic field with the smallest discriminant over which $E_{10}$ has rational points which are not cusps and therefore $Q(\sqrt{-1})$ is the quadratic field over which there exists a quartic polynomial whose square root has a continued fraction of period $9$.

Furthermore $\rank(E_{10}(Q(\sqrt{2})))=1$ whereas $\rank(E_{10}(K))=0$ for any quadratic field $K$ with $|\Delta|<8$. In fact the rank of $E_{10}(\Q(\sqrt{-2}))$ is $0$. Thus $Q(\sqrt{2})$ is the quadratic field with the smallest discriminant over which there are infinitely many quartic polynomials $f(x)$ such that the period of the continued fraction expansion of $\sqrt{f(x)}$ is 9.
\end{Proof}

\begin{Theorem}
\label{thm:period11}
The quadratic field $\Q\left(\sqrt{-1}\right)$ is the quadratic field with the smallest $|\Delta|$ over which there is a square free quartic polynomial $f(x)$ for which the continued fraction of $\sqrt{f(x)}$ has period $11$, whereas $\Q\left(\sqrt{5}\right)$ is the quadratic field with the smallest $|\Delta|$ over which there are infinitely many such quartic polynomials.
\end{Theorem}
\begin{Proof}
Let $f(x)$ be a quartic polynomial. If the point $\infty^+-\infty^-$ has order $12$ on the curve $y^2=f(x)$ then the period of the continued fraction of $\sqrt{f(x)}$ is either $11$ or $22$. An elliptic curve with a torsion point of order $12$ has the following parametrization $E_{b,c}:y^2+(1-c)xy-by=x^3-bx^2$ where
\begin{eqnarray*}
(b,c)=\left(\frac{t(2t-1)(2t^2 -2t+1)(3t^2 -3t+ 1)}{(t-1)^4},-\frac{t(2t-1)(3t^2 -3t+ 1)}{(t-1)^3}\right),
\end{eqnarray*}
where $t(t-1)(2t-1)(2t^2-2t+1)(3t^2-3t+1)\ne 0$.
A quartic model for the elliptic curve $E_{b,c}$ on which $\infty^+-\infty^-$ is a point of order $12$ is given by $y^2=(X^2+u_{12})^2-4v_{12}(X+w_{12})$ where
\begin{eqnarray*}
u_{12}(t)=\frac{12t^8 -120t^7 +336t^6 -468t^5 +372t^4 -168t^3 + 36t^2 -1}{4(t-1)^6},\\
v_{12}(t)=-\frac{t(2t-1)(2t^2 -2t+1)(3t^2 -3t+ 1)}{(t-1)^4},\; w_{12}=\frac{6t^4 -8t^3 + 2t^2 + 2t-1}{2(t-1)^3},
\end{eqnarray*}
with $\displaystyle k_{12}(t)=\frac{4t(2t-1)^2(3t^2 -3t+ 1)^3}{(t-1)^{11}}$. Now the continued fraction of $y$ can be found in \cite{ellipticfunctionfields} up to $a_{11}(y)$:
{\scriptsize\begin{eqnarray*}
y=\Big[X^2 +u_{12}(t), \frac{(t-1)^4}{2t(2t-1)(2t^2-2t+1)(3t^2-3t+1)}\left(X - \frac{6t^4-8t^3+2t^2+2t-1}{2(t-1)^3}\right)
,2\left(X + \frac{6t^4-10t^3+8t^2-4t+1}{2(t-1)^3}\right)
,\\ - \frac{(t-1)^3}{2t(2t-1)(3t^2-3t+1)}\left(X - \frac{2t^4+2t^3-6t^2+4t-1}{2(t-1)^3}\right)
,- \frac{2(3t^2−3t+1)}{(t-1)(2t^2-2t+1)}\left(X + \frac{2t^4-4t^3+6t^2-4t+1}{2(t-1)^3}\right)
,\\\frac{(t-1)^6}{2t(2t-1)(3t^2-3t+1)^2}\left(X − \frac{(2t-1)(2t^2-2t+1)}{2(t-1)^3}\right)
,\frac{2(2t-1)(3t^2-3t+1)}{(t-1)^5}\left(X - \frac{(2t-1)(2t^2-2t+1)}{2(t-1)^3}\right)
,\\-\frac{(t-1)^{10}}{2t(2t-1)^2(2t^2-2t+1)(3t^2-3t+1)^2}\left(X + \frac{2t^4-4t^3+6t^2-4t+1}{2(t-1)^3}\right)
,\frac{-2(2t-1)(3t^2-3t+1)^2}{(t-1)^8}\left(X - \frac{2t^4+2t^3-6t^2+4t-1}{2(t-1)^3}\right)
,\\ \frac{2(t-1)^{11}}{2t(2t-1)^2(3t^2-3t+1)^3}\left(X + \frac{6t^4-10t^3+8t^2-4t+1}{2(t-1)^3}\right)
,\frac{2(2t-1)(3t^2-3t+1)^2}{(t-1)^7(2t^2-2t+1)}\left(X - \frac{6t^4-8t^3+2t^2+2t-1}{2(t-1)^3}\right)
,\\\frac{(t-1)^{11}}{2t(2t-1)^2(3t^2-3t+1)^3}\left(X^2 +u_{12}(t)\right),\ldots
\Big].
\end{eqnarray*}}
One concludes that for $\mu y$ to have a periodic continued fraction of period $11$ one must have $a_{11}(\mu y)/a_0(\mu y)=2$, in other words $k_{12}(t)=1/\mu^2$, see Remark \ref{rem1}. This implies that $1/\mu$ is the $y$-coordinate of a rational point on the curve $\displaystyle y^2=\frac{4t(2t-1)^2(3t^2 -3t+ 1)^3}{(t-1)^{11}}$. Considering the square free part of $k_{12}(t)$, one needs to find a rational point on the elliptic curve $C_{12}:y^2=t(3t^2-3t+1)(t-1)$. The Jacobian elliptic curve of $C_{12}$ is $E_{12}:y^2=x^3+x^2+x.$

One has $E_{12}(\Q(\sqrt{-3}))=\Z/2\Z\times\Z/2\Z$, where all of the points are cusps. In fact $E_{12}(\Q(\sqrt{-1}))=\Z/8\Z$, where at least one of the points is not a cusp. Therefore $\Q(\sqrt{-1})$ is the quadratic field with the smallest absolute value of its discriminant $|\Delta|$ over which there is a quartic polynomial for which the continued fraction of the square root is periodic with period $11$.
 One has $\rank(E_{12}(\Q\sqrt{d}))=0$ where $d=-1,\pm2,\pm3,-5$, whereas $\rank(E_{12}(\Q(\sqrt{5})))=1.$ Therefore the quadratic field with the smallest $|\Delta|$ over which $E_{12}$ has a rational point is $Q(\sqrt{5})$. It follows that $\Q(\sqrt{5})$ is the quadratic field with the smallest $|\Delta|$ over which there are infinitely many square free quartic polynomials for which the continued fraction of the square root is periodic with period $11$.
\end{Proof}
\begin{Remark}
In \cite{ellipticfunctionfields} the periods 18 and 22 were shown to be realized over $\Q$. More precisely the authors used the parametrization of elliptic curves with torsion points of order $10$ and $12$ to write the square free quartic polynomials $f(x)$ over $\Q$ such that the continued fraction expansion of $\sqrt{f(x)}$ is of period 18 and 22 respectively.
\end{Remark}
 Now one has the following direct consequence giving a complete list of periods occurring over $\Q(\sqrt{-1})$ and $\Q(\sqrt{-3})$.
\begin{Corollary}
 Assume that $f(x)$ is a square free quartic polynomial over $K$ whose leading coefficient is a square. If the continued fraction expansion of $\sqrt{f(x)}$ is periodic with period $n$ then
 \begin{myitemize}
 \item[i.] If $K=\Q$ or $\Q(\sqrt{-3})$ then \[n\in S=\{1,2,3,4,5,6,7,8,10,14,18,22\}.\]
 Moreover for any $r\in S$ there is a square free quartic polynomial defined over $\Q$ such that the continued fraction of $\sqrt{f(x)}$ is periodic with period $r$.
 \item[ii.] If $K=\Q(\sqrt{-1})$ then \[n\in T= S\cup\{9,11\}.\]
 Moreover for any $r\in T$ there is a square free quartic polynomial defined over $\Q(\sqrt{-1})$ such that the continued fraction of $\sqrt{f(x)}$ is periodic with period $r$.
  \end{myitemize}
\end{Corollary}
\begin{Proof}
The list $S$ in (i) was proved to be complete over $\Q$, see the discussion in \S 3. The classification of torsion points of an elliptic curve defined over $\Q(\sqrt{-1})$ and $\Q(\sqrt{-3})$ is given in Theorem \ref{thm:Qi}. More specifically a torsion point on $E$ has maximal order $12$. Therefore the possible periods occurring over $\Q(\sqrt{-1})$ and $\Q(\sqrt{-3})$ are those occurring over $\Q$ together with $9,11$. According to Theorem \ref{thm:period9}, the field with the smallest absolute value of its discriminant over which the periods 9 and 11 occur is $\Q(\sqrt{-1})$.
\end{Proof}

\section{The even periods 10, 12, 14}

 The parametrization of elliptic curves with torsion points over quadratic fields was performed explicitly in \cite{Rabarison}. The parameters are rational points on a modular curve whose genus is either 1 or 2.

We write explicitly the square free quartic polynomials $f(x)$ such that the continued fraction of $\sqrt{f(x)}$ is periodic of even period $n$ where $n$ occurs over a quadratic field but not over the rational field. In this section we find those that give rise to the periods $ 10, 12$ and $14$. Although the periods $10$ and $14$ occur over $\Q$, there are different families of quartic polynomials over quadratic fields that give rise to the periods $10$ and $14$.

\subsection{Period 10}

Over $\Q$, quartic polynomials whose square roots have periodic continued fraction with period $10$ originate from
elliptic curves with torsion points of order $6$. Over quadratic fields, there is another source for these quartic polynomials, namely, elliptic curves whose torsion is of order $11$.
\begin{Theorem}
\label{thm:period10}
Let $f(x)$ be the quartic polynomial $f(x)=(x^2+u_{10})^2-4v_{10}(x+w_{10})$ defined over a quadratic field $K$ where
\begin{myenumerate}
\item $u_{10}(t)=(3t^2+6t-1)/4,\;v_{10}(t)=4t(t+1),\;w_{10}(t)=-(t-1)/2,$ where $t(t+1)(3t^2+6t-1)\ne0$, \\
or;
\item \begin{eqnarray*} u_{10}(s,t)&=& -\frac{s^4 + 2 s^3 t + t^2 + 6 s t^2 - 3 s^2 t (2 + t)}{4 t^2};\\v_{10}(s,t)&=&-\frac{(-1+s)s(s-t)}{t};\\w_{10}(s,t)&=&-\frac{s^2-t-st}{2t}\end{eqnarray*}
    where $(t,s)\in X_1(11)(K)$, $X_1(11):s^2-s=t^3-t^2$,
and, $t(t-1)(t^5-18t^4+35t^3-16t^2-2t+1)\ne 0$.
\end{myenumerate}
Then the continued fraction of $\sqrt{f(x)}$ is periodic with period $10$.
\end{Theorem}
\begin{Proof} In order for the continued fraction expansion of $\sqrt{f(x)}$ to be periodic of period $10$, the order of the curve $y^2=f(x)$ at $\infty^+-\infty^-$ must be either $6$ or $11$, see Proposition \ref{prop:periodlength}.

i. The curve $C:y^2=f(x)$ has order $6$ at $\infty^+-\infty^-$, this can be found for example in \cite[p. 109]{ellipticfunctionfields}. ii. An elliptic curve with a torsion point of order $11$ can be described as follows
\begin{eqnarray*}
y^2+\frac{(st+t-s^2)}{t}xy+\frac{s(s-1)(s-t)}{t}y=x^3+\frac{s(s-1)(s-t)}{t}x^2
\end{eqnarray*}
where $(t,s)\in X_1(11):s^2-s=t^3-t^2$,
and, $t(t-1)(t^5-18t^4+35t^3-16t^2-2t+1)\ne 0$. Observing that $\displaystyle (b,c)=\left(-\frac{s(s-1)(s-t)}{t},1-\frac{st+t-s^2}{t}\right)$, the result follows from Proposition \ref{prop:parametrizing quartics}.
\end{Proof}

\subsection{Period 12}

\begin{Theorem}
\label{period12}
Let $f(x)$ be the quartic polynomial $f(x)=(x^2+u_{12})^2-4v_{12}(x+w_{12})$ defined over a quadratic field $K$ where {\scriptsize\begin{eqnarray*}u_{12}(t,s)&=&-\frac{1}{16t^{10}}\Big\{1 - 2 t - 5 t^2 + 12 t^3 + 23 t^4 - 56 t^5 - 21 t^6 + 84 t^7 +
  12 t^8 - 58 t^9 + 11 t^{10} + 16 t^{11} - 2 t^{12} - 4 t^{13} + t^{14} \\&+&
  s^2 (-1 + t)^4 (-1 + t + t^2)^2 -
  2 s (-1 + t)^2 (1 - 2 t - 3 t^2 + 11 t^3 - 6 t^4 + 2 t^5 - 3 t^6 -
     6 t^7 - t^8 + t^9)\Big\};\\v_{12}(t,s)&=&-\frac{
 1 - 6 t + 14 t^2 - 13 t^3 + t^4 + 2 t^5 + 4 t^6 - t^7 - t^8 +
    s (-1 + t)^2 (-1 + 4 t - 5 t^2 + 2 t^4 + t^5)}{2t^9};\\w_{12}(t,s)&=&-\frac{-1 + 5 t - 9 t^2 + 5 t^3 + 2 t^4 - 3 t^5 - 2 t^6 + t^7 -
 s (-1 + t)^2 (-1 + t + t^2)}{4 t^5} \end{eqnarray*}}
  where $(t,s)$ is a point on $X_1(13)(K),\;X_1(13):s^2=t^6-2t^5+t^4-2t^3+6t^2-4t+1$ such that $t(t-1)(t^3-4t^2+t+1)\ne0$. Then the continued fraction of $\sqrt{f(x)}$ is periodic with period $12$.
  Further the quadratic field with the smallest $|\Delta|$ over which such quartic polynomial exists is $\Q(\sqrt{17})$. Moreover there are at most finitely many such quartic polynomials over a quadratic field $K$.
\end{Theorem}
\begin{Proof}
If the curve $C:y^2=f(x)$ has a point of order $13$ at $\infty^+-\infty^-$ then the continued fraction of $\sqrt{f(x)}$ is periodic with period $12$. An elliptic curve with a point of order $13$ is one of the following elliptic curves $y^2+(1-c)xy-by=x^3-bx^2$ where $(0,0)$ is a point of order 13 and
\begin{eqnarray*}
1-c&=&\frac{(t-1)^2(t^2+t-1)s-t^7+2t^6+3t^5-2t^4-5t^3+9t^2-5t+1}{2t^5};\\
b&=&-\frac{(t-1)^2\left((t^5+2t^4-5t^2+4t-1)s-t^8-t^7+4t^6+2t^5+t^4-13t^3+14t^2-6t+1\right)}{2t^9},
\end{eqnarray*}
where $(t,s)$ is a point on $X_1(13):s^2=t^6-2t^5+t^4-2t^3+6t^2-4t+1$ such that $t(t-1)(t^3-4t^2+t+1)\ne0$. The field $\Q(\sqrt{17})$ is the quadratic field with the smallest $|\Delta|$ such that $X_1(13)(K)$ contains a $K$-rational point which is not a cusp, see \cite[Theorem 3]{Najman}. Moreover the finiteness follows because $X_1(13)$ is a curve of genus $2$.
\end{Proof}
\begin{Remark}
Over $\Q(\sqrt{17})$, one has $(2,\sqrt{17})\in X_1(13)\left(\Q(\sqrt{17})\right)$. Therefore the continued fraction of the square root of the following quartic polynomial is periodic over $\Q(\sqrt{17})$ with period $12$
\[f(x)=\left(x^2-3 \frac{219 + 221 \sqrt{17}}{8192}\right)^2-\frac{107 - 51 \sqrt{17}}{256}\left(x+\frac{ 51 + 5 \sqrt{17}}{128}\right).\]
\end{Remark}

\subsection{Period 14}

\begin{Theorem}
\label{thm:period14}
Let $f(x)$ be the quartic polynomial $f(x)=(x^2+u_{14})^2-4v_{14}(x+w_{14})$ defined over a quadratic field $K$ where
\begin{myenumerate}
\item {\scriptsize\[u_{14}(t)=-4 - \frac{1}{4 t^2} + \frac{2}{t} + t + t^2,\;v_{14}(t)=  (-1 + t) (-1 + 2 t),\;w_{14}(t)=2- \frac{1}{2 t} - t;\]}
or;
\item {\scriptsize\begin{eqnarray*}u_{14}(t,s)&=&\frac{-1}{4 (1 + t)^6 (1 + t + t^2)^2}\Big\{1 + (8 - 6 s) t + (30 - 14 s + s^2) t^2 -
  2 (-39 + 11 s + s^2) t^3 + (143 - 16 s + s^2) t^4 \\&+& 4 (49 + s) t^5+(199 + 12 s) t^6 + 2 (70 + 3 s) t^7 + 63 t^8 + 14 t^9 + t^{10} \Big\};\\v_{14}(t,s)&=&\frac{ -t^3 (1 + t) (1 + t + 3 t^2 + t^3) -
    s t (-1 - t - 2 t^2 + t^4)}{(1 + t)^6 (1 + t + t^2)};\\w_{14}(t,s)&=&\frac{1 - (-4 + s) t + (7 + s) t^2 + 9 t^3 + 5 t^4 + t^5}{
 2 (1 + t)^3 (1 + t + t^2)}\end{eqnarray*}}
 where $(t,s)\in X_1(15)(K),\;X_1(15):s^2 + ts + s = t^3 + t^2$,
and $t(t+1)(t^2+t+1)(t^4+3t^3+4t^2+2t+1)(t^4-7t^3-6t^2+2t+1)\ne 0$.
\end{myenumerate}
Then the continued fraction of $\sqrt{f(x)}$ is periodic with period $14$.
\end{Theorem}
\begin{Proof}
For the continued fraction expansion of $\sqrt{f(x)}$ to be periodic of period $14$, the order of the curve $y^2=f(x)$ at $\infty^+-\infty^-$ must be either $8$ or $15$, see Proposition \ref{prop:periodlength}. In what follows we describe both possibilities respectively.

i. An elliptic curve has order $8$ at $\infty^+-\infty^-$ if it is parametrized as follows \[y^2+\left(1-(t-1)(2t-1)/t\right)xy-(t-1)(2t-1)y=x^3-(t-1)(2t-1)x^2.\]  ii. An elliptic curve with a torsion point of order $15$ can be described as follows $y^2+(1-c)xy-by=x^3-bx^2$ where
\begin{eqnarray*}
1-c&=&\frac{(t^2-t)s+(t^5+5t^4+9t^3+7t^2+4t+1)}{(t+1)^3(t^2+t+1)};\\
b&=&-\frac{t(t^4-2t^2-t-1)s+t^3(t+1)(t^3+3t^2+t+1)}{(t+1)^6(t^2+t+1)}
\end{eqnarray*}
where $(t,s)\in X_1(15):s^2 + ts + s = t^3 + t^2$,
and, $t(t+1)(t^2+t+1)(t^4+3t^3+4t^2+2t+1)(t^4-7t^3-6t^2+2t+1)\ne 0$.

In (i) and (ii) one uses Proposition \ref{prop:parametrizing quartics} to produce the square free quartic polynomial $f(x)$.
\end{Proof}

\begin{Remark}
In Theorem \ref{thm:period10} and Theorem \ref{thm:period14} the smallest field over which the periods 10 and 14 occur is the rational field $\Q$.
\end{Remark}

\section{The periods 13, 26, 15, 30, 17, 34}

In this section we parametrize the square free quartic polynomials $f(x)$ defined over a quadratic field $K$ for which the continued fraction expansion of $\sqrt{f(x)}$ is periodic with period 26, 30 or 34. Furthermore we present necessary and sufficient conditions in order for the continued fraction expansion of $\mu\sqrt{f(x)},\;\mu\in K,$ to be periodic with period 13, 15 or 17. The reason we investigate the pairs $(13, 26), (15, 30), (17, 34)$ of periods is that each pair arises from elliptic curves of the same torsion order. For example elliptic curves with $K$-rational torsion points of order $14$ give rise to quartic polynomials of either periods 13 or 26, see Proposition \ref{prop:periodlength}.

\subsection{Periods 13 and 26}
\begin{Theorem}
\label{period26}
Let $f(x)=(x^2+u)^2-4v(x+w)$ be defined over a quadratic field $K$ where
{\scriptsize\begin{eqnarray*}u(t,s)&=&\frac{-1}{4 (1 - 3 t^2 - t^3 + t^4)^2}(1 - 12 t^2 + s^2 (-1 + t)^4 t^2 +
   8 t^3 + 14 t^4 - 4 t^5 - 4 t^7 + t^8 -
   2 s t (3 - 8 t + t^2 + 12 t^3 - 7 t^4 - 2 t^5 + t^6));\\v(t,s)&=&\frac{(-1 + t) t (s (-1 + 2 t - 2 t^3) + t (-1 + t - t^3 + t^4))}{(1 -
  3 t^2 - t^3 + t^4)^2};\\w(t,s)&=&\frac{1 - s (-1 + t)^2 t - 4 t^2 + t^4}{2 (1 - 3 t^2 - t^3 + t^4)},\end{eqnarray*}}
   $(t,s)\in X_1(14)(K),\; X_1(14):s^2+st+s=t^3-t$ such that $t(t-1)(t+1)(t^3-9t^2-t+1)(t^3-2t^2-t+1)\ne0$. Then one has:
  \begin{myitemize}
  \item[i.] The continued fraction of $\sqrt{f(x)}$ is periodic with period $26$. The field $\Q(\sqrt{-7})$ is the quadratic field with the smallest $|\Delta|$ over which there are quartic polynomials whose square root has a periodic continued fraction with period $26$.
  \item[ii.] There exists a $\mu\in K$ such that the continued fraction expansion of $\mu\sqrt{f(x)}$ is of period $13$ if and only if $\alpha_{13}(t,s)$ is a $K$-square for some $(t,s)\in X_1(14)(K)$, $t(t-1)(t+1)(t^3-9t^2-t+1)(t^3-2t^2-t+1)\ne0$, where
  {\scriptsize\begin{eqnarray*}
  \alpha_{13}(t,s)&=&(-1 + t) t(1 + t) (1 - t - 2 t^2 + t^3)  (-s - t + s t) (-s - s^2 - t + s t + 3 s^2 t + 2 t^2 +
   4 s t^2 - 3 s^2 t^2 + 3 t^3 - 4 s t^3 + s^2 t^3 \\&-& 2 t^4 +
   s t^4) (s + t - 2 s t - t^2 + 2 s t^3 + t^4 - t^5) (-s^2 - s^3 -
   2 s t + 3 s^2 t + 5 s^3 t - t^2 + 8 s t^2 - 9 s^3 t^2 + 4 t^3 -
   10 s t^3 \\&-& 11 s^2 t^3 + 5 s^3 t^3 - 5 t^4 - 5 s t^4 + 18 s^2 t^4 +
   4 s^3 t^4 - t^5 + 35 s t^5 - 5 s^2 t^5 - 6 s^3 t^5 + 21 t^6 -
   18 s t^6 - 8 s^2 t^6 + 2 s^3 t^6 - 7 t^7 \\&-& 14 s t^7+ 10 s^2 t^7 -
   13 t^8 + 17 s t^8 - 5 s^2 t^8 + 6 t^9 - 5 s t^9 + s^2 t^9).
   \end{eqnarray*} }
   \end{myitemize}
\end{Theorem}
\begin{Proof}
An elliptic curve with a torsion point of order $14$ is given by $y^2+(1-c)xy-by=x^3-bx^2$ where
\begin{eqnarray*}
1-c&=&\frac{t^4-st^3+(2s-4)t^2-st+1}{(t+1)(t^3-2t^2-t+1)};\\
b&=& -\frac{-t^7+2t^6+(2s-1)t^5+(-2s-1)t^4+(-2s+2)t^3+(3s-1)t^2-st}{(t+1)^2(t^3-2t^2-t+1)^2},
\end{eqnarray*}
where $(t,s)\in X_1(14):s^2+st+s=t^3-t$ such that $t(t-1)(t+1)(t^3-9t^2-t+1)(t^3-2t^2-t+1)\ne0$. Now one uses Proposition \ref{prop:parametrizing quartics} to produce $f(x)$.

One has $X_1(14)(\Q(\sqrt{-7}))=\Z/2\Z\times\Z/6\Z$ where not all of these torsion points are cusps, moreover $X_1(14)$ contains only $K$-rational points which are cusps over any quadratic field $K$ of smaller $|\Delta|$, see for example \cite{Najman}. 

According to Lemma \ref{lem:periodicexpansion} and Remark \ref{rem1} one deduces that $k(t,s)$ in $\sqrt{f(x)}=\left[ a_0, \overline{a_1,\ldots,a_{r-1},2a_0/k(t,s),k(t,s)a_1,\ldots,a_{r-1}/k(t,s),2a_0}\right]$ is a $K$-square if and only if $\mu \sqrt{f(x)}$ has period $13$ where $\mu^2=k(t,s)$. One can obtain $k(t,s)$ using Lemma \ref{lem:periodicexpansion}, for example using Remark \ref{rem:partial quotients}, $\displaystyle k(t,s)=\frac{a_7}{a_6}=\frac{b_6(x-c_7)}{b_7(x-c_6)}$ where $c_6=c_7$. Now $\displaystyle k(t,s)=\frac{v(t,s)\beta(t,s)^2}{s_7}$ for some polynomial $\beta(t,s)$. Therefore $k(t,s)$ is a square if and only if $v(t,s)/s_7$ is a square. Using the definition of $s_7$ in Remark \ref{rem:partial quotients} and clearing the denominator of $v(t,s)/s_7$ one takes $\alpha_{13}(t,s)$ to be the square free part. Hence $k(t,s)$ is a square if and only if $\alpha_{13}(t,s)$ is a square.

\end{Proof}

\subsection{Periods 15 and 30}
\begin{Theorem}
\label{period30}
Let $f(x)=(x^2+u)^2-4v(x+w)$ be defined over a quadratic field $K$ where
 {\scriptsize\begin{eqnarray*}u(t,s)&=&\frac{-1}{4 (1 + t)^{10} (-1 + 2 t + t^2)^2}\Big\{1 + 34 t - 69 t^2 - 60 t^3 - 3 t^4 - 570 t^5 + 943 t^6 + 792 t^7 +
 1315 t^8 + 1262 t^9 + 273 t^{10} + 36 t^{11} \\&+& 95 t^{12} + 42 t^{13} + 5 t^{14}+ s^2 (1 - 7 t + 6 t^2 - 2 t^3 + t^4 + t^5)^2 +2 s (-3 + 22 t + 6 t^2 - 126 t^3 + 151 t^4 - 276 t^5 + 116 t^6 +
    84 t^7 \\&-& 53 t^8 + 30 t^9 + 38 t^{10} + 10 t^{11} + t^{12}))\Big\};\\v(t,s)&=&-\frac{ (-1 + t)^3 (s (-1 + 8 t - 2 t^2 + 8 t^3 + 3 t^4) +
    t (5 - 8 t + 15 t^2 + 11 t^4 + 8 t^5 + t^6))}{(1 + t)^8 (-1 +
    2 t + t^2)};\\w(t,s)&=&\frac{-1 - 7 t + 13 t^2 + 7 t^3 + 33 t^4 + 15 t^5 + 3 t^6 + t^7 +
  s (1 - 7 t + 6 t^2 - 2 t^3 + t^4 + t^5)}{
 2 (1 + t)^5 (-1 + 2 t + t^2)} \end{eqnarray*}} where $(t,s)\in X_1(16)(K),\;X_1(16):s^2=t(t^2+1)(t^2+2t-1)$ and $t(t-1)(t+1)(t^2+1)(t^2-2t-1)(t^2+2t-1)\ne0.$ Then one has:
 \begin{myitemize}
 \item[i.] The continued fraction of $\sqrt{f(x)}$ is periodic with period $30$. The field $\Q(\sqrt{-15})$ is the quadratic field with the smallest $|\Delta|$ over which there is a quartic polynomial whose square root has a periodic continued fraction with period $30$. Moreover there are at most finitely many such quartic polynomials over a quadratic field $K$.
  \item[ii.]  There exists a $\mu\in K$ such that the continued fraction expansion of $\mu\sqrt{f(x)}$ is of period $15$ if and only if $\alpha_{15}(t,s)$ is a $K$-square for some $(t,s)\in X_1(16)(K)$, $t(t-1)(t+1)(t^2+1)(t^2-2t-1)(t^2+2t-1)\ne0$, where
  {\scriptsize\begin{eqnarray*}
  \alpha_{15}(t,s)&=&(1 + t) (-1 + 2 t + t^2) (-s^3 - s t + 5 s^2 t + 17 s^3 t +
   5 t^2 + 20 s t^2 - 45 s^2 t^2 - 90 s^3 t^2 - 69 t^3 - 216 s t^3 -
   12 s^2 t^3 \\&+& 114 s^3 t^3 + 470 t^4 + 1384 s t^4 + 888 s^2 t^4 +
   153 s^3 t^4 - 1786 t^5 - 3844 s t^5 - 1447 s^2 t^5 + 111 s^3 t^5 +
   3598 t^6 + 4340 s t^6 \\&+& 723 s^2 t^6 + 276 s^3 t^6 - 4018 t^7 -
   4160 s t^7 - 1936 s^2 t^7 + 156 s^3 t^7 + 4430 t^8 + 4488 s t^8 -
   1144 s^2 t^8 + 129 s^3 t^8 - 3474 t^9 \\&+& 2090 s t^9 -
   1081 s^2 t^9 + 95 s^3 t^9 - 1048 t^{10} + 2428 s t^{10} -
   1175 s^2 t^{10} + 38 s^3 t^{10} - 1996 t^{11} + 3304 s t^{11} -
   548 s^2 t^{11}\\ &+& 18 s^3 t^{11} - 3102 t^{12} + 1304 s t^{12} -
   256 s^2 t^{12} + 7 s^3 t^{12} - 1006 t^{13} + 732 s t^{13} - 101 s^2 t^{13} +
    s^3 t^{13} - 526 t^{14} + 380 s t^{14}\\ &-& 15 s^2 t^{14} - 62 t^{15} +
   48 s t^{15} + 250 t^{16} - 8 s t^{16} + 122 t^{17} - s t^{17} + 19 t^{18} +
   t^{19}) (s^5 + s^3 t - 13 s^4 t - 32 s^5 t - 12 s^2 t^2 \\&+&
   35 s^3 t^2 + 393 s^4 t^2 + 421 s^5 t^2 + t^3 + 69 s t^3 +
   287 s^2 t^3 - 1241 s^3 t^3 - 4839 s^4 t^3 - 2920 s^5 t^3 -
   133 t^4 - 1841 s t^4 \\&-& 3007 s^2 t^4+ 14737 s^3 t^4 +
   31199 s^4 t^4 + 11415 s^5 t^4 + 2624 t^5 + 20091 s t^5 +
   12842 s^2 t^5 - 92859 s^3 t^5 - 113607 s^4 t^5 \\&-& 25416 s^5 t^5-
   23200 t^6 - 106863 s t^6 + 16038 s^2 t^6 + 350211 s^3 t^6 +
   241891 s^4 t^6 + 33947 s^5 t^6 + 105093 t^7 + 253022 s t^7 \\&-&
   356848 s^2 t^7 - 827395 s^3 t^7 - 328917 s^4 t^7 - 34336 s^5 t^7 -
   237961 t^8 - 54966 s t^8 + 1300808 s^2 t^8 + 1255195 s^3 t^8\\ &+&
   343061 s^4 t^8 + 23370 s^5 t^8 + 204632 t^9 - 854374 s t^9 -
   2229106 s^2 t^9 - 1368280 s^3 t^9 - 259986 s^4 t^9 -
   8160 s^5 t^9 \\&+& 107496 t^{10}+ 1766206 s t^{10} + 2726998 s^2 t^{10} +
   1243556 s^3 t^{10} + 86922 s^4 t^{10} - 430 s^5 t^{10} - 489355 t^{11} -
   2857553 s t^{11} \\&-& 3173649 s^2 t^{11} - 557066 s^3 t^{11} -
   8950 s^4 t^{11} + 5200 s^5 t^{11} + 1247799 t^{12} + 3839261 s t^{12} +
   1399993 s^2 t^{12} - 30086 s^3 t^{12}\\ &-& 65338 s^4 t^{12} -
   2930 s^5 t^{12} - 1686328 t^{13} - 1618311 s t^{13} - 350716 s^2 t^{13} -
   21942 s^3 t^{13} + 40482 s^4 t^{13} + 1040 s^5 t^{13} \\&+& 863640 t^{14} +
   1583195 s t^{14} + 603068 s^2 t^{14} - 626330 s^3 t^{14} -
   8202 s^4 t^{14} - 10 s^5 t^{14} - 1303543 t^{15} - 483020 s t^{15} \\&+&
   2536536 s^2 t^{15} - 332694 s^3 t^{15} + 8998 s^4 t^{15} -
   672 s^5 t^{15} - 139293 t^{16} - 3155332 s t^{16} + 2113176 s^2 t^{16} -
   392826 s^3 t^{16} \\&+& 17338 s^4 t^{16} - 91 s^5 t^{16} + 950128 t^{17} -
   3451268 s t^{17} + 2602492 s^2 t^{17} - 331599 s^3 t^{17} +
   6863 s^4 t^{17} - 192 s^5 t^{17} \\&+& 1419280 t^{18} - 5186604 s t^{18} +
   2217984 s^2 t^{18} - 174197 s^3 t^{18} + 6461 s^4 t^{18} -
   135 s^5 t^{18} + 2944963 t^{19} - 4980589 s t^{19} \\&+& 1438957 s^2 t^{19} -
   111725 s^3 t^{19} + 3949 s^4 t^{19} - 40 s^5 t^{19} + 3275313 t^{20} -
   3930407 s t^{20} + 931395 s^2 t^{20} - 58043 s^3 t^{20} \\&+&
   1419 s^4 t^{20} - 21 s^5 t^{20} + 3117520 t^{21} - 2882435 s t^{21} +
   483570 s^2 t^{21} - 22335 s^3 t^{21} + 629 s^4 t^{21} - 8 s^5 t^{21} +
   2645648 t^{22}\\ &-& 1693321 s t^{22} + 202366 s^2 t^{22} - 8889 s^3 t^{22} +
   215 s^4 t^{22} - s^5 t^{22} + 1811983 t^{23} - 838386 s t^{23} +
   78936 s^2 t^{23} - 2839 s^3 t^{23} \\&+& 31 s^4 t^{23} + 1076709 t^{24} -
   365798 s t^{24} + 24800 s^2 t^{24} - 497 s^3 t^{24} + s^4 t^{24} +
   553144 t^{25} - 126518 s t^{25} + 5206 s^2 t^{25}\\& -& 42 s^3 t^{25} +
   229512 t^{26} - 32626 s t^{26} + 726 s^2 t^{26} - 2 s^3 t^{26} +
   76551 t^{27} - 6439 s t^{27} + 69 s^2 t^{27} + 20445 t^{28} - 933 s t^{28} \\&+&
   3 s^2 t^{28} + 4008 t^{29} - 81 s t^{29} + 504 t^{30} - 3 s t^{30} +
   35 t^{31} + t^{32})(-1 + t) (-s + 5 t + 8 s t - 8 t^2 - 2 s t^2 + 15 t^3 + 8 s t^3\\ &+&
   3 s t^4 + 11 t^5 + 8 t^6 + t^7).
  \end{eqnarray*}}
 \end{myitemize}
\end{Theorem}
\begin{Proof}
An equation describing an elliptic curve with a torsion point of order $16$ is given by $y^2+(1-c)xy-by=x^3-bx^2$ where
\begin{eqnarray*}
1-c&=&\frac{(t-1)(t^4+2t^3+6t-1)}{(t+1)^5(t^2+2t-1)}s+\frac{t^5+t^4+14t^3+6t^2+9t+1}{(t+1)^5};\\
b&=&-\frac{(t-1)^3(3t^4+8t^3-2t^2+8t-1)}{(t+1)^8(t^2+2t-1)}s-\frac{t(t-1)^3(t^2+1)(t^4+8t^3+10t^2-8t+5)}{(t+1)^8(t^2+2t-1)}
\end{eqnarray*}
where $(t,s)\in X_1(16):s^2=t(t^2+1)(t^2+2t-1)$ and $t(t-1)(t+1)(t^2+1)(t^2-2t-1)(t^2+2t-1)\ne0.$ In order to produce $f(x)$ one needs to apply the transformations in Proposition \ref{prop:parametrizing quartics}.

That $ X_1(16)(\Q(\sqrt{-15}))$ is the quadratic field with the smallest $|\Delta|$ containing rational points which are not cusps can be found in \cite{Najman}. The finiteness result follows because $X_1(16)$ is of genus 2.

In view of Lemma \ref{lem:periodicexpansion} and Remark \ref{rem1} one deduces that $k(t,s)$ in the expansion $\sqrt{f(x)}=\left[ a_0, \overline{a_1,\ldots,a_{r-1},2a_0/k(t,s),k(t,s)a_1,\ldots,a_{r-1}/k(t,s),2a_0}\right]$ is a $K$-square if and only if $\mu \sqrt{f(x)}$ has period $15$ where $\mu^2=k(t,s)$. One can obtain $k(t,s)$ using Lemma \ref{lem:periodicexpansion}, for example using Remark \ref{rem:partial quotients}, $\displaystyle k(t,s)=\frac{a_7}{a_8}=\frac{b_8(x-c_7)}{b_7(x-c_8)}$ where $c_7=c_8$. Now $\displaystyle k(t,s)=b_8/b_7=\frac{\beta(t,s)^2}{s_8v(t,s)}$ for some polynomial $\beta(t,s)$. Using the definition of $s_8$ in Remark \ref{rem:partial quotients} and clearing the denominator of $\displaystyle \frac{1}{s_8v(t,s)}$ one defines $\alpha_{15}(t,s)$ to be the square free part. Now $k(t,s)$ is a $K$-square if and only if $\alpha_{15}(t,s)$ is a $K$-square.
\end{Proof}
\subsection{Periods 17 and 34}
\begin{Theorem}
\label{period34}
Let $f(x)=(x^2+u)^2-4v(x+w)$ where {\scriptsize\begin{eqnarray*}u(t,s)&=&\frac{-1}{16 t^6 (-1 - 3 t + t^3)^2}\Big\{(s^2 (2 + 3 t + t^2)^2 -
  6 s t (-2 - 11 t - 20 t^2 - 20 t^3 - 9 t^4 + 5 t^5 + 7 t^6 +
     2 t^7) \\&+&
  t (8 + 49 t + 150 t^2 + 239 t^3 + 254 t^4 + 167 t^5 + 6 t^6 -
     90 t^7 - 120 t^8 - 79 t^9 - 12 t^{10} + 4 t^{11}))\Big\};\\v(t,s)&=&\frac{ (1 + 2 t + 2 t^2 + t^3) (-1 - t - t^2 + 4 t^3 + 6 t^4 + 4 t^5 +
    3 t^6 + t^7 + s (1 - t + 2 t^3 + t^4))}{2t^5 (-1 - 3 t + t^3)^2};\\w(t,s)&=&-\frac{s (2 + 3 t + t^2) +
 t (5 + 11 t + 11 t^2 + 10 t^3 + t^4 - 2 t^5)}{4 t^3 (-1 - 3 t + t^3)}
   \end{eqnarray*}} where $(t,s)\in X_1(18)(K),\;X_1(18):s^2=t^6+2t^5+5t^4+10t^3+10t^2+4t+1$ and $t(t+1)(t^2+t+1)(t^3-3t-1)\ne 0$. Then one has
   \begin{myitemize}
   \item[i.] The field $\Q(\sqrt{33})$ is the quadratic field with the smallest $|\Delta|$ over which there are quartic polynomials whose square root has a periodic continued fraction with period $34$. Over a quadratic field $K$ there are at most finitely many such quartic polynomials.
  \item[ii.] There exists a $\mu\in K$ such that the continued fraction expansion of $\mu\sqrt{f(x)}$ is of period $17$ if and only if $\alpha_{17}(t,s)$ is a $K$-square for some $(t,s)\in X_1(18)(K)$, $t(t+1)(t^2+t+1)(t^3-3t-1)\ne 0$, where
   {\scriptsize\begin{eqnarray*}
  \alpha_{17}(t,s)&=&(-1 - 3 t + t^3)(1 + t) (2 s + 5 t + s t + 6 t^2 + 3 t^3 + t^4) (-1 + s - t - s t -
   2 t^2 + 3 s t^2 + 9 t^3 + 5 s t^3 + 21 t^4 + s t^4 + 14 t^5\\ &+&
   4 t^6 + t^7) (4 s^2 - 4 s^3 + 4 t + 12 s t - 4 s^2 t - 4 s^3 t +
   37 t^2 + 41 s t^2 - 17 s^2 t^2 - s^3 t^2 + 138 t^3 + 64 s t^3 -
   2 s^2 t^3 - 4 s^3 t^3\\ &+& 235 t^4 + 54 s t^4 - 37 s^2 t^4 -
   16 s^3 t^4 + 80 t^5 - 205 s t^5 - 174 s^2 t^5 - 25 s^3 t^5 -
   578 t^6 - 720 s t^6 - 239 s^2 t^6 - 19 s^3 t^6\\ &-& 1370 t^7 -
   801 s t^7 - 160 s^2 t^7 - 7 s^3 t^7 - 1183 t^8 - 234 s t^8 -
   69 s^2 t^8 - s^3 t^8 + 216 t^9 + 213 s t^9 - 24 s^2 t^9 +
   1424 t^{10} \\&+& 213 s t^{10}- 6 s^2 t^{10} + 1475 t^{11} + 100 s t^{11} -
   s^2 t^{11} + 912 t^{12} + 38 s t^{12} + 428 t^{13} + 9 s t^{13} + 157 t^{14} +
   s t^{14} + 41 t^{15}+ 8 t^{16} \\&+& t^{17}) (16 s^4 - 16 s^5 + 24 s^2 t +
   112 s^3 t - 72 s^4 t - 32 s^5 t + 4 t^2 + 108 s t^2 +
   468 s^2 t^2 + 172 s^3 t^2 - 296 s^4 t^2 + 8 s^5 t^2 + 166 t^3\\ &+&
   1240 s t^3 + 1988 s^2 t^3 - 304 s^3 t^3 + 70 s^4 t^3 +
   152 s^5 t^3 + 1483 t^4 + 5975 s t^4 + 3894 s^2 t^4 +
   1314 s^3 t^4 + 2071 s^4 t^4 + 239 s^5 t^4 \\&+& 6776 t^5 +
   15176 s t^5 + 9404 s^2 t^5 + 11520 s^3 t^5 + 4284 s^4 t^5 +
   80 s^5 t^5 + 17831 t^6 + 27326 s t^6 + 33578 s^2 t^6 +
   27812 s^3 t^6 \\&+& 3459 s^4 t^6 - 142 s^5 t^6 + 29373 t^7 +
   51974 s t^7 + 82930 s^2 t^7 + 33300 s^3 t^7 + 57 s^4 t^7 -
   178 s^5 t^7 + 35080 t^8 + 110326 s t^8\\ &+& 127366 s^2 t^8 +
   20106 s^3 t^8 - 2154 s^4 t^8 - 88 s^5 t^8 + 47247 t^9 +
   201895 s t^9 + 129842 s^2 t^9 + 2454 s^3 t^9 - 1909 s^4 t^9 -
   21 s^5 t^9 \\&+& 97308 t^{10} + 289240 s t^{10} + 95864 s^2 t^{10} -
   5290 s^3 t^{10} - 864 s^4 t^{10} - 2 s^5 t^{10} + 200855 t^{11} +
   332136 s t^{11} + 58224 s^2 t^{11} \\&-& 4360 s^3 t^{11} - 241 s^4 t^{11} +
   326410 t^{12} + 317611 s t^{12} + 33652 s^2 t^{12} - 1772 s^3 t^{12} -
   43 s^4 t^{12} + 411028 t^{13} + 257510 s t^{13}\\ &+& 19260 s^2 t^{13} -
   438 s^3 t^{13} - 4 s^4 t^{13} + 411057 t^{14} + 177292 s t^{14} +
   9976 s^2 t^{14} - 60 s^3 t^{14} + 334919 t^{15} + 103132 s t^{15} \\&+&
   4204 s^2 t^{15} - 2 s^3 t^{15} + 226785 t^{16} + 49927 s t^{16} +
   1334 s^2 t^{16} + 128639 t^{17} + 19583 s t^{17} + 302 s^2 t^{17} +
   60753 t^{18} \\&+& 6080 s t^{18} + 46 s^2 t^{18} + 23634 t^{19} +
   1472 s t^{19} + 4 s^2 t^{19} + 7527 t^{20} + 264 s t^{20} + 1933 t^{21} +
   31 s t^{21} + 382 t^{22} + 2 s t^{22}\\ &+& 55 t^{23} + 5 t^{24}).
  \end{eqnarray*}}
       \end{myitemize}
\end{Theorem}
\begin{Proof}
The Tate's equation describing an elliptic curve with a point of finite order $18$ is given by $y^2+(1-c)xy-by=x^3-bx^2$ where
{\scriptsize\begin{eqnarray*}
1-c&=&\frac{-t^2-3t-2}{2t^3(t^3-3t-1)}s-\frac{-2t^5+t^4+10t^3+11t^2+11t+5}{2t^2(t^3-3t-1)} ;\\
b&=&\frac{(t+1)(t^2+t+1)(t^4+2t^3-t+1)}{2t^5(t^3-3t-1)^2}s+\frac{(t+1)(t^2+t+1)(t^7+3t^6+4t^5+6t^4+4t^3-t^2-t-1)}{2t^5(t^3-3t-1)^2}
\end{eqnarray*}}
where $(t,s)\in X_1(18):s^2=t^6+2t^5+5t^4+10t^3+10t^2+4t+1$ and $t(t+1)(t^2+t+1)(t^3-3t-1)\ne 0$. One uses the transformations in Proposition \ref{prop:parametrizing quartics} to produce $f(x)$. That $ X_1(18)(\Q(\sqrt{33}))$ is the quadratic field with the smallest $|\Delta|$ containing rational points which are not cusps can be found in \cite{Najman}. The finiteness result follows because $X_1(18)$ is of genus 2.

 According to Lemma \ref{lem:periodicexpansion} and Remark \ref{rem1} one deduces that $k(t,s)$ in the expansion $\sqrt{f(x)}=\left[ a_0, \overline{a_1,\ldots,a_{r-1},2a_0/k(t,s),k(t,s)a_1,\ldots,a_{r-1}/k(t,s),2a_0}\right]$ is a $K$-square if and only if $\mu \sqrt{f(x)}$ has period $17$ where $\mu^2=k(t,s)$. One can obtain $k(t,s)$ using Lemma \ref{lem:periodicexpansion}, for example using Remark \ref{rem:partial quotients}, $\displaystyle k(t,s)=\frac{a_9}{a_8}=\frac{b_8(x-c_9)}{b_9(x-c_8)}$ where $c_8=c_9$. Now $\displaystyle k(t,s)=b_8/b_9=\frac{s_9\beta(t,s)^2}{v(t,s)}$ for some polynomial $\beta(t,s)$. Using the definition of $s_9$ in Remark \ref{rem:partial quotients} and clearing the denominator of $s_9/v(t,s)$ one sets $\alpha_{17}(t,s)$ to be the square free part. Now $k(t,s)$ is a $K$-square if and only if $\alpha_{17}(t,s)$ is a $K$-square.
\end{Proof}

\hskip-18pt\emph{\bf{Acknowledgements.}} All calculations were performed using \Sage and \Mathematica.

\bibliographystyle{plain}
\footnotesize
\bibliography{ContinuedFractions}
Department of Mathematics and Actuarial Science\\ American University in Cairo\\ mmsadek@aucegypt.edu
\end{document}